\documentclass[11pt]{amsart}
\usepackage[margin=1in]{geometry}

\usepackage{amssymb}
\usepackage{amsthm}
\usepackage{amsmath}
\usepackage{mathrsfs}
\usepackage{amsbsy}
\usepackage{bm}
\usepackage{hyperref}
\usepackage{tikz}
\usepackage{array}
\usepackage{enumerate}
\usepackage{enumitem}
\usepackage{bbm}
\usepackage{comment}
\usepackage{mathtools}
\usepackage{makecell} 
\usepackage{colortbl}
\usepackage{xcolor}
\usepackage{comment}

\DeclareFontFamily{U}{mathx}{}
\DeclareFontShape{U}{mathx}{m}{n}{<-> mathx10}{}
\DeclareSymbolFont{mathx}{U}{mathx}{m}{n}
\DeclareMathAccent{\widecheck}{0}{mathx}{"71}

\definecolor{LightBlue}{rgb}{0,0.8,1} 
\definecolor{lavender}{rgb}{0.4,0,1}

\hypersetup{colorlinks=true, citecolor=LightBlue, linkcolor=lavender,urlcolor=lavender}

\usepackage{hhline}
\allowdisplaybreaks
\usepackage[noadjust]{cite}

\usepackage{caption}
\usepackage[noabbrev,capitalise,nameinlink]{cleveref}
\crefname{conjecture}{Conjecture}{Conjectures}

\theoremstyle{definition}

\usepackage{etoolbox}

\definecolor{darkgreen}{RGB}{0, 120, 0}

\begin{document}

\title[All Quiet on the AI Front]{All Quiet on the AI Front: \\ The Future of Undergraduate Math Research, \\ Viewed from Duluth}
\date{}
\subjclass[2010]{}

\author[]{Colin Defant}
\address[]{Department of Mathematics, Harvard University, Cambridge, MA 02138, USA}
\email{colindefant@gmail.com} 

\author[]{Eliot Hodges}
\address[]{Fine Hall, Princeton University, Princeton, NJ 08540, USA}
\email{eliotljhodges@gmail.com} 

\author[]{Noah Kravitz}
\address[]{St John's College, Oxford and Mathematical Institute, University of Oxford; St Giles', Oxford OX1 3JP, UK}
\email{noah.kravitz@maths.ox.ac.uk}

\maketitle

\section{Backdrop}

In this time of uncertainty and anxiety about artificial intelligence (AI) in mathematics, it is more important than ever to reflect on how we can prepare the next generation of professional research mathematicians for success.  In this article, we report, from the front lines, on how AI is affecting the long-running summer math research program at the University of Minnesota Duluth. 
We do not have answers to the big questions about how AI will shape the future of professional mathematics. What we can offer instead is a collection of practical observations about the role of AI in the 2026 Duluth program: how it was useful and where it presented potential pitfalls; what surprised us and how we adapted; and what went well and what we wished we had done differently.

This article is aimed at educators, particularly those who supervise undergraduate research programs. Because the Duluth program resembles PhD-level research, much of our discussion also applies to graduate mentoring, with the caveat that the longer graduate timescale allows for more flexibility. Student readers may pick up pointers on how to make more effective use of AI tools in their research experiences.

Student research training, once the exclusive prerogative of graduate programs, is increasingly common for undergraduates. In standard course work, students learn about foundational mathematical concepts and build technical acuity. Undergraduate research experiences provide controlled environments in which students can begin developing the independence and initiative that they will later need to make new discoveries and set research agendas. Early research exposure also provides a primer on para-mathematical skills such as effective communication. Ideally, students come away from such experiences with a notion of how professional math research builds on and goes beyond classroom learning. 

The Duluth Research Experience for Undergraduates (REU) program, founded by Joe Gallian in 1977, has for decades been a premier venue for undergraduate research. 321 students have participated in the program to date, and many alumni now hold positions at top schools.  The program is designed to approximate a PhD research experience. Each student works individually on an open combinatorics (or combinatorics-adjacent) problem from the recent literature and presents progress in weekly board talks. Students are encouraged to consult one another and the program's research advisors, but each remains responsible for advancing their own project. In a typical year, almost every student produces at least one original research paper for publication in a well-regarded peer-reviewed journal. 

The dramatis personae of the 2026 REU program were: ten undergraduate students; two directors (Joe Gallian and Colin Defant); four research advisors (Eliot Hodges, Claire Kaneshiro, Noah Kravitz, and Carl Schildkraut); and twenty short-term visitors (mostly recent program alumni). Students and directors resided in Duluth for the duration of the 9-week program. Advisors stayed for 3-6 weeks each and facilitated both research and social aspects of the program.  To protect our students' privacy, we have anonymized examples and removed identifying details about mathematical projects.  The reader should bear in mind that what we frame as general pronouncements are directly informed by our summer 2026 observations.

\section{Goals and Guiding Considerations}

Our first priority as mentors was ensuring that AI would enhance the research experience rather than detract from it. This included teaching students to think critically about the pros and cons of choosing to use AI tools for various tasks. Going into 2026, we knew from previous summers that AI-free undergraduate research experiences are highly beneficial, so one possibility on the table was an indiscriminate program-wide ban on AI. We quickly discarded this option, however, for two reasons. The first, more ideological, reason was that AI is sure to play a major role in the future of math research, so preparing students to be professional mathematicians includes actively equipping them to use AI effectively. The second, more practical, reason was that an AI ban would have been nearly impossible to implement fairly and would have undercut trust-building in the mentor-mentee relationship.

Once we decided to allow AI use, we had to figure out which types of use to encourage. We were concerned that querying of AI might inadvertently abbreviate or bypass pedagogically valuable elements of the research process. We worried that even if AI played only an auxiliary role, students might become over-reliant on these tools and miss out on opportunities to develop their own problem-solving skills. 

To assuage these apprehensions, we emphasized to our students that being stuck is an important part of the research process. The friction of straining against a difficult part of a proof is what most often leads to insights. We suggested that students wait to consult AI on problem-solving matters until after they had made a reasonable effort on their own---for instance, one day's worth of solid thinking. Even though it may soon become standard practice to outsource ``routine'' lemmas to AI assistants, we emphasized that working out such details oneself is a good way to build technique. For less problem-solving-oriented tasks, such as literature search and generation of code and examples, we encouraged students to take full advantage of the speed-ups offered by AI. We reminded students that they needed to carefully check AI output (especially potentially spurious citations) and make sure that they understood everything they were using.

We initially expected that students would incline toward overuse of AI tools. To our surprise, we observed a wide variety of attitudes toward AI, ranging from antipathy to indifference to cautious adoption. Only a minority of students seemed to be frequent users. Many were unaware of the mathematical capabilities of the top AI models and needed instruction on harnessing the strength of their preferred AI tools (e.g., submitting queries on ``high thinking'' settings).

Finally, our students got caught up in the general AI-induced realignment of the mathematical community's value system. Now that AI models can mass-produce substantial proofs, many results that formerly might have been considered solid incremental progress are increasingly going the way of routine AI cleanup. The bar for publication of human work is accordingly shifting (mostly rising). Even if they were not fully aware of the nuances of this shift, many of our Duluth students expressed preoccupations with the perceived ``importance'' of their summer research results---noticeably more in the 2026 cohort than in previous years' groups. Just a few years ago, ambitious students commonly fixated on \emph{quantity} of output as a way to impress graduate admissions committees or win awards. Although this mindset was certainly counterproductive, we have now seen the pendulum swing too far in the other direction: several 2026 students felt that low-hanging fruit was, so to speak, not worth picking at all. 

We reminded our students that the primary goal of the research program is for students to gain familiarity with the research process as a whole, both the experience of prolonged effort on a difficult problem and the surrounding research soft skills. The perceived quality of research output has always been only a secondary consideration---a flashy result is a nice bonus if it happens, but fixating on outcomes is counterproductive. An early research experience is ideally a vehicle for acquiring long-term skills. We think that this pragmatic training-centric attitude towards undergraduate (and perhaps also graduate) research will necessarily become more common as AI capabilities continue improving and the bar for human novelty continues rising. 

We maintain that students who refrain from substantive AI use in their research can be just as successful as those who do partake. Many of our students in 2026 solved difficult problems despite making (at most) minimal use of AI tools. Even though AI systems are capable of dispatching many undergraduate-level problems, there is still plenty of room for human creativity and ingenuity. At the end of the summer, we were quite proud of our students' work and accomplishments.

\section{Gradations of Student AI Use}

When deployed carefully, AI tools can enhance the research process by allowing researchers to focus their efforts on the most salient aspects of the problems at hand. AI use becomes fraught, however, when it starts \emph{replacing} human thinking instead of \emph{supplementing} it. In general, researchers with more experience and mathematical maturity can afford to outsource higher-level tasks to the machine without jeopardizing their autonomy; the dangers of overuse are heightened for novices who do not yet have solid foundations for thinking independently. Where a seasoned expert might quickly zero in on a gap in an AI-generated proof, a novice could spend a week checking technical details before finding the crucial error. Likewise, it takes experience to identify the ``core idea'' underlying lengthy calculations output by an AI model or to discern whether a suggested future direction is worth pursuing. Based on our observations at the 2026 Duluth program, we have identified several paradigms for more and less helpful patterns of student AI use.

AI tools can effectively accelerate many of the auxiliary tasks of math research, such as searching the literature, writing code, and generating examples to stress-test conjectures. Doing these activities by hand of course has its benefits---you never know where you will stumble across an important observation---but time saved by AI assistance can also be valuable.

Explaining one's ideas to an interlocutor, even a non-expert, is a great way to take stock of progress, identify potential issues, and brainstorm next steps. In Duluth, these interlocutors have traditionally been advisors and visitors. Many students in 2026 supplemented this human feedback by using their AI assistants for ``rubber duck debugging.'' Even though AI assistants were typically less helpful than advisors and visitors, they had the advantage of being available 24/7 to listen. 

Like all research, undergraduate research often requires working with unfamiliar tools and new topics on the fly. Since AI systems have a broad mathematical background, they can help students fill gaps in their knowledge while they are finding their sea legs. Such assistance is particularly useful for streamlining background reading and extracting technical statements from thorny parts of the literature. A potential danger, however, is that reliance on AI for foundational knowledge can lead students to get in over their heads. This type of AI use can also make it harder for mentors to detect when a student needs extra support or a proposed project is too ambitious. 

After obtaining a mathematical result, one should always look for ways to simplify proofs and further strengthen theorems. AI is very effective at these kinds of extensional tasks. In 2026, several of our students used AI tools to streamline their arguments, and one used AI to find a ``matching example'' demonstrating the optimality of an inequality they proved. Such ex post facto substantive consultation of AI does not detract from the problem-solving experience. 

Substantive consultation of AI presents more of a mixed bag \emph{during} the problem-solving phase. On the one hand, using AI to handle technical lemmas can free an experienced researcher to focus on the big picture. On the other hand, for students such ``routine'' steps can be valuable sites of learning and acclimatization to unfamiliar techniques. We think that students stand to gain more from their research experiences if they minimize their use of AI for shortcutting proofs and instead get their hands dirty themselves.

More pernicious outcomes can arise in matters of judgment. Since one never wants to waste time pursuing dead ends, it can be tempting to use AI systems as quick ``temperature checks'' for gauging the difficulty of conjectures, the feasibility of proof strategies, or the interest of proposed lines of inquiry. In 2026, we saw some students put significant effort into chasing down unpromising conjectures suggested by AI, while others were led to concede defeat too quickly after AI systems deemed problems difficult. On the flip side, some students found success in places that AI systems had deemed unpromising; it is important to take subjective judgments from AI with big grains of salt. Developing a reliable sense for promising directions is just as important a skill as working out the details of the resulting problems.

The worst scenario imaginable is for a student to outsource their entire project to an AI system. Besides having lost the benefit of grappling with difficult concepts, such a student can end up with an uninteresting, unmotivated paper whose content they do not understand. Despite the pressure to produce a research paper with limited time in a research program, students will be better off in the long run if they prioritize the development of their own problem-solving skills. It goes without saying that forthrightness about AI use (at any level) allows mentors to provide students with the most effective support.

\section{Pedagogical Considerations}

\subsection{Problem Selection}

The biggest single decision in student research mentoring is problem selection. A student who likes their assigned problem is much more likely to find success! The ideal project has the potential for progress at a variety of levels---a paper with several different open questions and room for partial results provides more flexibility than a paper with only a single all-or-nothing conjecture. At the Duluth program, we aim for easily motivated problems with low barriers to entry so that students can get to problem-solving quickly, without spending too much of the summer on background reading. We typically source most of our problems by trawling the arXiv over the course of the preceding year.

For the 2026 program, we committed to assigning only problems that were not easily solvable by AI. A month in advance, we ran some of our selected papers through ChatGPT 5.4 Pro and asked it to try to solve the stated open problems. With minimal prompting (e.g., ``keep trying'' a few times), ChatGPT settled eight of the twenty problems that we gave it. Many of these solutions were nontrivial and would have made for excellent undergraduate-level (or even PhD-level) papers if they had been discovered by humans. Faced with this assault on our usual source of projects, we considered several possible responses, each with its pros and cons.

The first response, and the one that we mostly adopted, was to continue as close as possible to business as usual. After throwing out the potential projects on which ChatGPT had made substantial progress, we still had enough remaining for our students. The problems that we ended up using were, on the whole, more difficult than usual, and for some topics we were scraping the bottom of the barrel earlier in the summer than we would have liked. The result was that, despite some student frustration with slow progress, we were mostly able to operate as in previous years. This approach will likely become less feasible in the future when stronger AI models can one-shot more of the problems that comprise the typical Duluth fare.

A related response was to pivot to seemingly more AI-proof projects in areas where AI models are not as capable. However, given the pace at which frontier models are improving, this is at best a short-term solution. A more promising pivot is turning to projects where the challenge is not only deploying technical tools to resolve preexisting conjectures but also creatively identifying new directions of inquiry. Even if an AI model later assists in answering such questions, the process of raising and sifting through them fosters independent inquiry and develops mathematical ``taste.'' Examples include finding tractable special cases of an out-of-reach conjecture and exploring properties of a recently defined object. In the 2026 Duluth program we used several such problems, especially of the latter variety, with good overall success. One student's topic ended up being so rich in fruitful directions that we eventually had to declare a moratorium on new results in the interest of getting to the writing stage before the end of the program. We caution, however, that such open-ended inquiry seems to work best with students of an appropriate temperament, and in 2026 
these projects worked less well for students who preferred more concrete goals. We also found that open-ended projects required extra preliminary vetting, since it was hard to evaluate quickly whether or not there were directions worth exploring.

A third possible response was to revisit our initial opposition to giving out problems that AI could solve. In Duluth we have hitherto taken the position that assigning genuinely open problems makes for a more authentic research experience and builds more independence. It is common in other settings, however, for mentors to assign undergraduate (and even graduate) research projects that they already more or less know how to solve. From a pedagogical perspective, it does not particularly matter whether a domain expert \emph{could} in theory resolve a student's assigned problem, and perhaps going forward we will have to take the same view of AI solutions. Of course, then one would have to grapple with a slew of further practical considerations, such as publishability of results and protecting students from getting accidentally ``scooped'' by AI assistants midway through the research process. We also mention the possibility that, rather than killing a project, an AI solution could serve as a springboard for a student's further investigations; we wish in hindsight that in 2026 we had not been so quick to discard problems on which ChatGPT had made progress.

Finally, one could embrace the coming storm and deliberately select projects that incorporate AI.  We hope to explore such ``alternative'' research styles in future summers, even though they might deviate substantially from the traditional models of undergraduate research.  We will have to keep an open mind about the details, however, because the volatility of the AI-and-math research landscape makes it difficult to make predictions for even a few months in the future.

\subsection{Writing and Editing}

The research process does not end when the main theorems are proven: one must then share one's findings with the broader mathematical community. New researchers often find the demands of professional mathematical writing to be among the biggest differences between solving a challenging homework assignment and producing a research paper.

At the Duluth program, we encourage all students to write up their results for publication in peer-reviewed journals. Each student initially writes and edits a complete paper draft on their own. Advisors and visitors (typically one of each) then read the draft and give detailed feedback on both mathematical content and writing style. Most students have never before been responsible for solo-authored research articles, and over the course of the multi-stage editing process they learn the ``dos and don'ts'' of writing a math paper.  We encourage students to think of revision as an opportunity to refine writing skills, not a tedious task to be checked off.

Providing such extensive feedback is time-intensive, and the full back-and-forth process can span several months. In 2026, we experimented with running student papers through a round of AI editing in advance of human editing. An AI proofreader that catches typos and flags grammatical mistakes can save human editors time and energy. Likewise, since an AI system can digest mathematical arguments in a matter of seconds, it can function as a first line of defense to check for soundness and preempt a human editor's request for additional explanation.

The second author (Eliot Hodges) used Claude to build a custom AI editing workflow. We initially trained the AI agent on a library of well-constructed past student papers in order to give it a feel for our preferred Duluth ``house style.'' The custom AI workflow takes a student paper's \LaTeX\ source code as input and then outputs a series of four text files with feedback at various ``scales.'' The first feedback file collects low-level mechanical and typographical issues, including typos, notational collisions, and grammatical mistakes. It also checks for recurrent writing concerns from prior years of the Duluth program, such as parallel structure and punctuation around centered equations. The second file records the results of a careful mathematical audit of the correctness of the paper and makes note of incorrect or otherwise suspicious steps. The third offers suggestions for how to improve readability and presentation on a local level, and the fourth gives broad-brushstroke advice about the paper's content and overall organization. 

One thing that our AI editing agent does \emph{not} do is implement fixes or otherwise execute rewrites. Whereas automatic correction would abet the continuation of bad habits, a list of comments creates an opportunity for a student to actively decide how (if at all) to address each proposed edit. This system also gives students room to preserve their individual writing styles. We believe that active incorporation of external feedback reinforces learning.

Although the new AI editing workflow added a step to our editing process, it quickly paid dividends. The drafts that human editors received after AI editing were noticeably more polished and easier to read. Since common pitfalls had already been mostly addressed, human editors could focus their energy on finer points of style and presentation. In at least one instance, the AI editor identified a gap in a proof that otherwise might have taken a human reader several weeks to discover. 

Using AI for editing of course also comes with pitfalls. The primary peril is that students who outsource editing to the machine can subsequently neglect careful revision of their own, in terms of both checking mathematical arguments and refining exposition. Students who place too much trust in the AI's editing judgment can feel tempted to produce sloppy writing and let the AI clean it up. In 2026, we observed some students ``edit'' their papers by asking an AI to flag errors, making the suggested changes, and repeating until the AI declared the work to be error-free. Such mindless editing undercuts authorial independence and bypasses the reflection demanded by ordinary writing and revision. Another shortcoming of AI editors is that LLMs are often better than humans at understanding poorly written or poorly motivated arguments. What an AI reader deems intelligible may not be clear for a human reader, so deferral to AI judgment on presentation can sometimes hurt human readability. 

One potential solution to student over-reliance on AI editors is setting explicit expectations that students edit their drafts to the best of their ability before turning to AI, just as in the past we asked students to polish their work before turning to human readers. 

\subsection{Training AI Literacy} 
We have mentioned the importance of experience and judgment for effectively interfacing with AI tools. Even though developing such mathematical maturity takes time, there are also concrete practical skills that yield immediate results. We believe that such ``AI literacy'' will be increasingly important in the future. We made some efforts to teach AI literacy in the 2026 Duluth program, but we wish we had emphasized it more.

We held an informal seminar on practical tips for optimizing AI output in mathematical contexts. Topics included turning on extended thinking modes, working in a local context (e.g., Codex or Claude Cowork), and rudimentary prompt engineering. Students were noticeably more confident AI users after the seminar. In addition, visitors Evan Chen and Janabel Xia ran a workshop on the basics of Lean and formal proof verification.

We also asked students to keep track of their interactions with AI throughout the summer, both so that they could determine which uses were most productive and so that they could accurately attribute credit when they later wrote papers. Finally, on a lighter note, we made up a quiz titled ``Duluthian vs. DuluthAIn'' in which students competed to distinguish between real abstracts from recent Duluth program papers and spoof abstracts generated by ChatGPT. This fun activity primed students to think critically about their engagement with AI-generated content in the literature.

\section{The Next Frontier}
Our experiences in 2026 have prepared us to help students make their AI use thoughtful and productive throughout the research process. We conclude this article with several tangible changes that we plan to make in the 2027 Duluth program. We believe that it is more important than ever for mathematicians to cultivate independence and mathematical taste.

As we discussed above, we will expand the range of projects that we assign. In addition to proposing more traditional projects with specific conjectures and sets of tools, we will keep an eye out for open-ended projects that call for student initiative in setting directions of inquiry.\footnote{The quantity and nature of arXiv postings are currently changing so dramatically that we may soon be forced to reconsider our strategy of systematically screening all postings in ``Combinatorics'' and related categories. Perhaps an AI assistant could help with the initial winnowing of potential project topics!} In the past, we have assigned projects primarily on the basis of the personal statements from students’ applications. Before the 2027 program begins, we will contact students individually to learn about their attitudes toward AI and the kinds of problems they would most enjoy. A more granular picture of each student's mathematical interests and temperament will help us select a style of project that is a good personal fit.

Advisors and visitors have always played an important role as sounding boards for students to discuss their work in progress. We hope that these more senior program members can also use their prior research experience to provide guidance on AI-related matters. For example, we will encourage students to share not only how their reasoning developed but also the role of AI in generating their ideas. Moreover, advisors and visitors can help students make sure that they fully understand content originating from AI tools. Similarly, when AI suggests a conjecture or a new research direction, visitors and advisors can help students assess its plausibility and identify fruitful alternatives. These conversations will give students valuable practice in vetting AI-generated ideas and exercising mathematical judgment.

An enduring goal of the Duluth program has been equipping students with the soft skills required for professional research mathematics. Examples in the past have included organizing material clearly during board presentations and formatting papers in \LaTeX. As it becomes increasingly necessary for working mathematicians to have AI-related skills, too, we will accordingly expand our student training. Starting in the 2027 Duluth program, we will introduce a weekly group workshop (perhaps dubbed \emph{Superintelligence Saturday}) devoted to developing students' AI literacy skills.

These plans are necessarily provisional because AI is sure to advance in the coming year. Our broader objective, however, will remain to help students take advantage of AI without sacrificing the independence, persistence, and judgment that a research experience is meant to cultivate.

\section*{Acknowledgments} 
We are indebted to Joe Gallian, who built the Duluth program and the wonderful community of Duluthians that we now cherish. We are also grateful to all of the members of the 2026 program for allowing us to experiment with, observe, and report on their interactions with AI. This includes the other two program advisors, Claire Kaneshiro and Carl Schildkraut. It also includes the 20 visitors: Glenn Bruda, Ruben Carpenter, Evan Chen, Mitchell Lee, Derek Liu, Annabel Ma, Emily Ma, Andrei Mandelshtam, Delong Meng, David Moulton, Jacob Paltrowitz, Isaac Rajagopal, Michael Ren, Pitchayut (Mark) Saengrungkongka, Maya Sankar, Robert Schneider, Kerry Seekamp, Natasha Ter-Saakov, Janabel Xia, and Sophie Zhu. 

C. Defant was supported in part by an AMS--Simons travel grant. E. Hodges was supported in part by the NSF Graduate Fellowship Research Program under grant DGE-2444107.  N. Kravitz was supported in part by a NSF Mathematical Sciences Postdoctoral Research Fellowship under grant DMS-2501336. Any opinions, findings, and conclusions or recommendations expressed in this material are those of the author(s) and do not necessarily reflect the views of the National Science Foundation. The 2026 Duluth REU was supported by Jane Street Capital, Axiom Math, and a donation from Larry Penn.   

\section*{Statement on Potential Conflicts of Interest} 
C. Defant started a residency at the AI-for-math company Axiom Math midway through the work leading to the present article.  His contributions are based solely on his observations at the Duluth REU and his experiences as an academic researcher; in particular, he has no intention to promote any specific AI tools or companies. 

Axiom Math also supported the 2026 Duluth program with an unrestricted donation that did not factor into the program's operations (including AI use) or the production of this article.

\end{document}